\documentclass[12pt]{article}

\usepackage{amssymb}

\def\C{{\mathbb C}}
\def\Aut{{\rm Aut}}

\usepackage{theorem}

\newtheorem{theorem}{Theorem}
\newtheorem{proposition}{Proposition}[section]

{\theorembodyfont{\rmfamily}
\newtheorem{definition}[proposition]{Definition}

}

\include{epsf}

\title{Enumeration of almost polynomial rational functions
with given critical values}

\author{Dmitri Panov\thanks{IPDE, IHES, Le Bois-Marie, 35, 
route de Chartres F-91440, Bures-sur-Yvette, France. 
E-mail: panov@ihes.fr\,.}, 
Dimitri Zvonkine\thanks{
Institut Math{\'e}matique de Jussieu,
Universit{\'e} Paris~VI, 175, rue du Chevaleret, 75013 Paris,
France.
E-mail: zvonkine@math.jussieu.fr\,.\/\/
The second author was partially suported by
EAGER - European Algebraic Geometry Research Training Network, 
contract No. HPRN-CT-2000-00099 (BBW) and by the Russian Foundation
of Basic Research grant 02-01-22004.}}

\date{\today}

\begin{document}

\maketitle

\begin{abstract}
Enumerating ramified coverings of the sphere with fixed
ramification types is a well-known problem first considered
by A.~Hurwitz~\cite{Hurwitz}. Up to now, explicit solutions
have been obtained only for some families of ramified coverings,
for instant, those realized by polynomials in one complex 
variable. In this paper we obtain an explicit answer for a large new
family of coverings, namely, the coverings realized by
simple almost polynomials, defined below. Unlike most other
results in the field, our formula is obtained by elementary methods.
\end{abstract}

\section{Rational functions and minimal factorizations of
permutations}

Let $f: \C \rightarrow \C$ be a rational function
of degree $n$ in one complex variable.

A {\em critical point} of $f$ is a point $z \in \C$ such that
$f'(z)=0$. Its {\em degree} is the number $a \geq 2$ such that $f$ 
looks like $f(z) = z^a$ in the neighborhood of the critical point.
A {\em critical value} of $f$ is its value at
a critical point. (Note that we do not count poles
as critical points.)

\begin{definition}
A rational function $f$ is called {\em simple} if every
critical value of $f$ has exactly one critical preimage. It is
called {\em almost polynomial} if the sum of orders of its
poles is smaller than the degree of each
critical point.
\end{definition}

Thus an almost polynomial rational function has a numerator
of a much bigger degree than the denominator.

Our goal is to find the number of simple almost polynomial
rational functions with fixed orders of poles and
fixed critical values of fixed
degrees, a problem first considered by V.~Arnold~\cite{Arnold}.

\bigskip

The general problem of enumerating ramified coverings of the sphere
with fixed ramification types can, in some sense, be solved
using the representation theory of the symmetric group; however
the answer one obtains is a rather complicated sum over the
irreducible representations. In purticular, there is still no
simple criterion allowing one to determine whether the number of coverings
is equal to $0$ or not.

For (not necessarily simple)
polynomials the problem can be reduced to a combinatorial problem
solved by I.~P.~Goulden, D.~M.~Jackson
in~\cite{GouJac}. Later, when the relation to polynomials was
discovered, their formula was reproved in
\cite{LanZvo}, \cite{PanZvo}, and~\cite{knizhka}
(appendix by D.~Zagier) by three different methods.
In~\cite{PanZvo} we also treated the case of almost polynomial functions
with just one simple pole. In all these cases, as well as in
the case considered in the present note, one obtains simple closed answers. 

A.~Goupil and G.~Schaeffer~\cite{GouSch} generalized 
Goulden and Jackson's result on polynomials to meromorphic
functions with a unique pole on Riemann surfaces of any genus,
but their answer is not as explicit.
For other results on the enumeration of ramified coverings
and their relation to the intersection theory on moduli
spaces see~\cite{ELSV}, \cite{Zvonkine} and the references
therein.

\bigskip

All rational functions we consider in the sequel are simple and
almost polynomial.

A rational function $f$ can be viewed as a ramified covering
of the Riemann sphere by the Riemann sphere. Going around
a critical value in the image we obtain a permutation of
the sheets in the preimage (the monodromy of the covering). 
It follows from Riemann's existence theorem 
that the problem of counting rational
functions can be reformulated in terms of permutations
(for details see~\cite{knizhka}, chapter~2). 

Let $\sigma \in S_n$ be a permutation of $n = \deg f$ elements.

\begin{definition}
A product $\sigma_k \dots \sigma_1 = \sigma$ is called a
{\em minimal factorization} of $\sigma$ if (i)~the group
generated by $\sigma_1, \dots, \sigma_k$ acts transitively
on the set $\{ 1, \dots, n \}$ and (ii)~the total number 
of cycles in the permutations 
$\sigma_1, \dots, \sigma_k, \sigma$ equals $kn-n+2$.
\end{definition}

Here $\sigma$ corresponds to the monodromy of $f$ at
$\infty$, while the $\sigma_j$'s are the monodromies
at the critical values. The conditions guarantee that the
ramified covering determined by the permutations
$\sigma_1, \dots, \sigma_k, \sigma$ is (i)~connected and
(ii)~of genus~$0$ (by the Riemann-Hurwitz formula). 
Rather than counting the cycles of the permutations,
it is more natural to consider their {\em defects}:
a defect being $n-\mbox{(the number of cycles)}$.
If the Euler characteristic of the covering surface
equals $\chi$, the sum of defects of the corresponding
monodromies equals $2n - \chi$. Hence
the total defect $2n-2$ of a minimal factorization is 
the smallest possible for a transitive factorization, which explains 
the word ``minimal''.

For shortness, we will say that {\em a 
rational function $f$ is of type} $(a_1, \dots, a_k)$,
$(p_1, \dots, p_c)$ if it has $c$ poles of orders
$p_1, \dots, p_c$ and $k$ critical points of degrees
$a_1, \dots, a_k$. (The sum $p=p_1 + \dots + p_c$ satisfies 
$p < a_j$ for all~$j$.) Similarly, we say that {\em a 
minimal factorization is of type} $(a_1, \dots, a_k)$,
$(p_1, \dots, p_c)$ if $\sigma$ has cycles of lengths
$p_1, \dots, p_c, n-p$, while each $\sigma_j$ is an $a_j$-cycle,
i.e., has exactly one cycle of length $a_j$ the other points
being fixed.

To a minimal factorization of type $(a_1, \dots, a_k)$,
$(p_1, \dots, p_c)$ we can assign a colored graph with oriented
edges called a {\em constellation}. It is obtained in the following
way. Take $n$ numbered vertices. For each $j$, $1 \leq j \leq k$
form an oriented polygon using the vertices from the cycle
of $\sigma_j$. The edges of the polygon are colored in  ``color'' $j$.
Now forget the numbers of the vertices. 

It is clear that the
constellation allows one to reconstitute the permutations
$\sigma_1, \dots, \sigma_k, \sigma$ up to a common conjugation.

\begin{definition} \label{Def:constellation}
A {\em constellation} is a connected graph whose edges are oriented and
colored in colors from $1$ to $k$, obtained by gluing together
$k$ oriented polygons with colors $1, \dots, k$ at some
of their vertices. A vertex of a polygon cannot be glued
to another vertex of the same polygon.
\end{definition}

To sum up, we now have three equivalent problems: given
$(a_1, \dots, a_k)$, $(p_1, \dots, p_c)$, count
the number of rational functions,
or of minimal factorizations, or of constellations
of this type. The three numbers differ by simple combinatorial
factors.

\bigskip

We will need some more remarks on constellations.

A constellation coming from a minimal
factorization has a natural embedding into a sphere
(see~\cite{knizhka}). It is
given by the following conditions: (i)~The orientations of the
edges determine the counterclockwise orientation on
each polygon; (ii)~If we choose any vertex and enumerate the 
colors of the polygons surrounding it in the counterclockwise order
starting from the smallest color we obtain an increasing 
sequence of colors.

\begin{figure}[h]
\begin{center}
\
\epsfbox{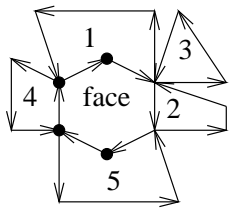}

\caption{The color of each polygon is marked inside it.
The essential vertices of the central face are shown in black;
these vertices form a cycle of $\sigma_5 \dots \sigma_1$.}
\label{Fig:faces}

\end{center}
\end{figure}

Cutting the sphere along the edges of the embedded constellation
we obtain $k+c+1$ pieces homeomorphic to open discs. Among
them, $k$ correspond to the polygons and the $c+1$ others
to the cycles of $\sigma$. The piece corresponding to the
long cycle (of length $n-p$) will be called the {\em
exterior face}; the other $c$ pieces will be called
{\em interior faces} or just {\em faces}. Note that
not all vertices of a face are contained in the corresponding
cycle of $\sigma$. Those which are will be called the
{\em essential vertices}, see Figure~\ref{Fig:faces}.

From now on we will always assume that the constellations
are embedded into a sphere. 

Note that the somewhat bizarre condition $p < a_j$ has a
simple interpretation for constellations: it means that
every polygon must have at least one edge on the exterior face.
However the algebraic meaning of this condition remains
mysterious.

\section{The main theorem}

Let $p_1, \dots, p_c$ be $c \geq 0$ positive integers
(orders of poles), $\sum p_i = p$.
Fix an integer $n > p$ (degree of the almost polynomial). 
Let $a_1, \dots a_k$ be $k \geq 2$ more 
positive integers (multiplicities of the critical points) 
satisfying $a_j >p$ for all $j$ and $\sum a_j = n+k+c-1$ 
(the Riemann-Hurwitz formula). Denote by $|\Aut \{ p_1, \dots, p_c \}|$
the number of permutations $s$ of $c$ elements such that
$p_i = p_{s(i)}$ for all $i$. For instance, 
$|\Aut \{ 4,4,3,2,2,2,2,2,1,1 \}| = 2! \cdot 1! \cdot 5! \cdot 2 !\,$.

Consider a permutation $\sigma \in S_n$ with cycle type
$(p_1, \dots, p_c, n-p)$.

\begin{theorem} \label{Thm:main}
The number of minimal factorizations of $\sigma$ into
$a_j$-cycles equals
$$
\frac{(k+c-2)!}{(k-2)!} \; p_1^2 \dots p_c^2 \; (n-p)^{k-1}. 
$$
The number of constellations as well as the number
of rational functions with fixed critical
values of type $(a_1, \dots, a_k)$,
$(p_1, \dots, p_c)$ equals
$$
\frac{1}{|\Aut \{ p_1, \dots, p_c \}|} \; 
\frac{(k+c-2)!}{(k-2)!} \; p_1 \dots p_c \; (n-p)^{k-2}. 
$$
\end{theorem}

The three assertions of the theorem are equivalent. Indeed,
the number of constellations and the number
of rational functions of a given type coincide, which follows from
Riemann's existence theorem (see~\cite{knizhka} for more
details). On the other hand,
to obtain a minimal factorization of $\sigma$ from
a given constellation, we must number
the vertices of the constellation in such a way that the product
$\sigma_k \dots \sigma_1$ equals $\sigma$. It is easy to 
see that there are
$$
p_1 \dots p_c (n-p) \cdot |\Aut \{ p_1, \dots p_c \}|
$$
ways to do that. (This number is also the number of permutations
that commute with $\sigma$.) 

In the next two sections we prove the theorem for constellations.

\section{Assembling a constellation}

We are going to prove the assertion of the theorem on the number
of constellations. We start by labeling the faces of
the constellations so as to make them distinguishable, which
kills the $|\Aut|$ factor.
We must show that the number of such constellations with
labeled faces equals
$$
\frac{(k+c-2)!}{(k-2)!} \; p_1 \dots p_c \; (n-p)^{k-2}.
$$

We proceed by induction on the number $c$ of
faces in the constellations. 

For $c=0$, the constellations have
no (interior) faces, i.e., they are ``trees'' glued of
a given set of polygons. Such constellations (also called
``cacti'') were enumerated in~\cite{GouJac},
\cite{LanZvo}, and~\cite{PanZvo}.
The answer one obtains is $n^{k-2}$.

\begin{figure}[h]
\begin{center}
\
\epsfbox{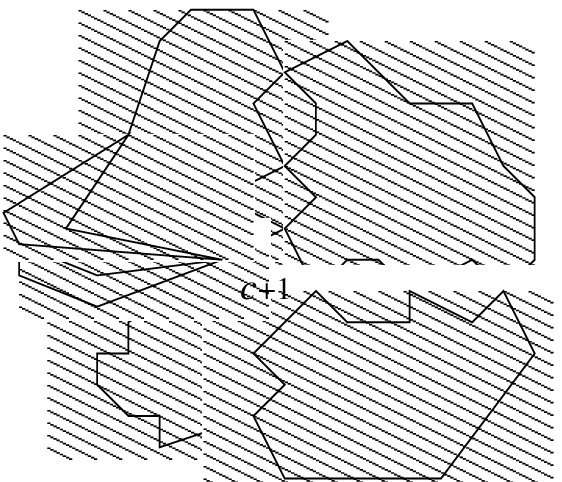}

\caption{The $(c+1)$st face in this figure has $4$
neighboring polygons and $5$ neighboring faces.}
\label{Fig:neighbors}
\end{center}
\end{figure}

Suppose the formula is established for constellations with
$\leq c$ faces and let us add one more face. 
A {\em polygon} of the constellation is a {\em neighbor} of
the $(c+1)$st face if it has at least one edge in
common with this face. A {\em face} of the constellation
is a {\em neighbor} of the $(c+1)$st face if its bounding
polygons are all neighbors of this face.
The $(c+1)$st face can have any number $2 \leq m \leq k$ of
neighboring polygons and any number $0 \leq d \leq c$ of
neighboring faces as in Figure~\ref{Fig:neighbors}.
(Note that a polygon that has only one vertex in common
with the face is not considered a neighbor.) The
condition $p < a_j$ implies that each neighboring face
is bounded by exactly two polygons.

To construct a constellation with $c+1$ faces we must
make the following choices.

{\bf 1.}~Choose the numbers $m$ and $d$.

{\bf 2.}~Choose $m$ polygons among $k$ and $d$ faces among $c$
to be the neighbors of the $(c+1)$st face. This gives
$$
{k \choose m} {c \choose d}
$$
choices. Denote by $D \subset \{1,\dots, c \}$ the set of the
$d$ neighboring faces.

\begin{figure}[h]
\begin{center}
\
\epsfbox{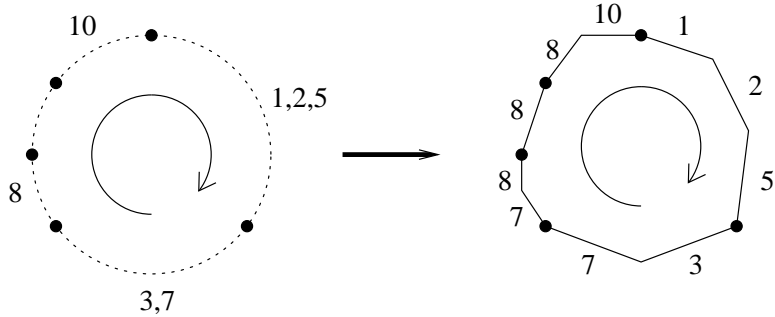}

\caption{How to obtain a face from a list of colors assigned
to each interval $V_i V_{i+1}$.}
\label{Fig:formface}
\end{center}
\end{figure}

{\bf 3.}~Form the $(c+1)$st face using the $m$ chosen polygons. This
is done in the following way. Denote the essential
vertices of the face by $V_1, \dots, V_{p_{c+1}}$.  
We must describe the colors and the order of the
edges that will form the intervals between $V_i$ and
$V_{i+1}$ for each $i$. We claim that such a disposition
of edges is uniquely determined once we have (arbitrarily)
assigned to each of the $m$ polygons the interval
$V_i V_{i+1}$ where its first edge will appear (as we go around the face in
the clockwise direction). Indeed, the disposition of
the edges can be obtained as follows. (a)~For each polygon
assigned to the interval $V_i V_{i+1}$ take one edge of
its color. (b)~Order these edges in the increasing order
of colors. (c)~In the case if either the biggest color used in
$V_{i-1} V_i$ is smaller than the smallest color in the
list for $V_i V_{i+1}$ or if there are no polygons assigned to
$V_i V_{i+1}$, add, at the beginning of the list of colors
for $V_i V_{i+1}$, the last color used in $V_{i-1}V_i$. 
An example of this algorithm is shown in Figure~\ref{Fig:formface}.

It is easy to see that the vertices $V_i$ are indeed the
essential vertices and that this is the unique way to achieve
this.

Thus the number of ways to form the $(c+1)$st face using $m$ given
polygons is $p_{c+1}^m$.

{\bf 4.}~Now we must choose the positions of the $d$ faces
that were chosen to be the neighbors of the $(c+1)$st face.
There are $m$ ``clefts''  between the polygons where these faces
can be placed. Several faces can appear in the same cleft.
In this case they must be ordered (starting from the face
closest to the $(c+1)$st one). Thus there are 
$$
\frac{(d+m-1)!}{(m-1)!}
$$
ways to choose the positions of the faces.

{\bf 5.}~Each of the neighboring faces has exactly two bounding
polygons. We now choose how many edges of the face will belong to
the polygon, say, on its left. For the $i$th face there are
$p_i$ choices. Thus we obtain a factor
$$
\prod_{i \in D} p_i .
$$
This product, of course, depends on the particular choice
of the $d$ neighboring faces. However we will soon see
that the remaining part of the formula contains the
complementary product
$$
\prod_{i \not\in D} p_i ,
$$
and thus it is not necessary to take a sum over all
possible choices.

{\bf 6.}~We have assembled together all the
polygons and faces that are neighbors of the $(c+1)$st face. They 
form a subconstellation $K$ of the total constellation. Now
we are going to {\em consider $K$ as a
unique polygon}. Indeed, we are going to show that the
remaining polygons are attached to $K$
in the same way as they would be attached to a unique polygon.
In particular, the number of ways to attach them is the same.
This will allow us to procede by induction.

Denote by $r = p_{c+1}+\sum_{i \in D} p_i$ the sum of lengths
of the faces of $K$ and by $v$ the total number of
vertices of $K$. These $v$ vertices are acted upon by $m$ cyclic
permutation corresponding to the polygons. The product
of these $m$ permutations splits the vertices into $d+2$ cycles
of lengths $p_i$, $i \in D$, $p_{c+1}$, and $v-r$.
The vertices of the last cycle will be called the 
{\em essential exterior vertices} of $K$.

Now consider a constellation $R$ formed by the $k-m$
polygons not used in $K$ and, in addition, one
$(v-r)$-gon of color $0$. Suppose that the faces of $R$ have
lengths $p_i$, $i \not\in D$. We are going to replace
the $(v-r)$-gon in $R$ by the constellation $K$.

Let
$p' = p_1 + \dots + p_c$ and $p = p_1 + \dots + p_{c+1}$.

{\bf 7.}~First we establish a one-to-one correspondence between
the vertices of the $(v-r)$-gon and the essential exterior
vertices of $K$ preserving their cyclic order. There
are $v-r$ ways of doing that, which will account for
a factor $v-r$ in the final formula. 

In principle, the number $v-r$ is different
for different choices of the $m$ polygons and the $d$
faces. However we will soon see that this number 
appears in the total sum only as a linear factor.
Therefore it will be a posteriori justified to replace
it by its average over the possible choices of $m$ polygons
and $d$ faces. 

The average value $\left< v-r \right>$ of $v-r$ equals
$$
\frac{m}{k} (n+k+c) -m - d -p_{c+1} - \frac{d}{c} p'
= \frac{mn+mc-k p_{c+1}}{k} -d \, \frac{p'+c}{c}.
$$
Indeed, the total number of edges in the $k$ polygons
equals $n+k+c$ by the Euler formula. Choosing randomly $m$ of
the $k$ polygons we obtain an average of $\frac{m}{k}(n+k+c)$
edges in the subconstellation $K$. Since $K$
has $d+1$ faces and $m$ polygons, its average number
of vertices is $\left< v \right> = \frac{m}{k}(n+k+c) - m-d$,
again by the Euler formula.
Now, $K$ has $d+1$ faces, one of which is always of
length $p_{c+1}$ and the other $d$ are chosen randomly
from $c$ possibilities. Thus $ \left< r \right> 
= p_{c+1} + \frac{d}{c} p'$, whence we obtain $\left< v-r \right>$.

{\bf 8.}~Let us go around the exterior face
of the constellation $K$ in the counterclockwise direction.
The colors of the edges we meet will be increasing 
in each interval between two consecutive essential exterior 
vertices. Then, as we pass an essential vertex, the number of the
color jumps down.

Suppose we are given a new polygon of some color $j$ that
does not appear in the constellation $K$. We want to
attach this polygon to the exterior of $K$ in such a way
that the cyclic order of colors around each vertex remains
increasing (see the remarks after 
Definition~\ref{Def:constellation}). It is easy to
see that in each interval between two consecutive
essential vertices $V_i$ (included) and $V_{i+1}$ (excluded)
there is a unique vertex to which the polygon can be attached.

Now, the essential exterior vertices of $K$ are in
a one-to-one correspondence with the vertices of a
$(v-r)$-gon in the new constellation $R$. We want to
replace the $(v-r)$-gon in $R$ by the subconstellation $K$.
To do that, for each polygon $j$ attached to some vertex
of the $(v-r)$-gon, we take the corresponding essential
exterior vertex $V_i$ of $K$ and attach our polygon to
the unique possible vertex between $V_i$ and $V_{i+1}$.

It is easy to see that the faces of $R$, even those that
have been modified by our operation, will still have the
same number of essential vertices as before. Thus there
is a unique way to substitute the $(v-r)$-gon in $R$
by the subconstellation $K$. The result is the constellation
we were trying to assemble.

{\bf 9.}~It remains to choose the constellation $R$.
By the induction assumption, there are
$$
\frac{(k-m-1+c-d)!}{(k-m-1)!} (n-p)^{k-m-1}
 \prod_{i \not\in D} p_i
$$
choices. Indeed, the constellation $R$ has $k-m$ polygons,
$c-d$ (interior) faces of lengths $p_i$, $i \not\in D$, 
and the length of its exterior face is
$n-p$ (the same as in the constellation that we
are assembling).

\section{Computing the sum}

The result of our investigation is that
the number of constellations (with labeled faces)
is given by the following sum:

$$
S=
p_1 \dots p_c
\sum_{m=2}^k {k \choose m} (n-p)^{k-m-1} p_{c+1}^{m-1}
\times \qquad \qquad  \qquad \qquad  \qquad \qquad
$$
$$
\times \sum_{d=0}^c {c \choose d}
\frac{(d+m-1)!}{(m-1)!}
\frac{(k-m-1+c-d)!}{(k-m-1)!}
\left(
\frac{mn \! + \! mc \! - \! k p_{c+1}}{k} -d \, \frac{p'+c}{c} 
\right).
$$

The rest of the proof is a sequence of elementary but rather
cumbersome computations.

We first compute the subsum over $d$ using the elementary relations
$$
\sum_{d=0}^c {a+d \choose d}{b-d \choose c-d} = {a+b+1 \choose c},
$$
$$
\sum_{d=0}^c d {a+d \choose d}{b-d \choose c-d} = 
(a+1) {a+b+1 \choose c-1}
$$
with $a = m-1$, $b = k+c-m-1$.
We find that the subsum is equal to
$$
\frac{(k+c-1)!}{(k-1)!} 
\left( \frac{mn+mc-k p_{c+1}}{k} \right)
-m \frac{(k+c-1)!}{k!} (p'+c)
$$
$$
= \frac{(k+c-1)!}{k!} [ m(n-p') - k p_{c+1}].
$$

Substituting this into the initial sum we obtain
$$
S = p_1 \dots p_c \frac{(k+c-1)!}{k!}
\sum_{m=2}^k {k \choose m} (n-p)^{k-m-1} p_{c+1}^{m-1}
[m(n-p') - k p_{c+1}].
$$
This sum can now be evaluated using two more elementary
identities
$$
\sum_{m=2}^k m {k \choose m} x^{k-m-1} y^{m-1} 
= k \left[ \frac{(x+y)^{k-1}}{x} - x^{k-2} \right]
$$
$$
\sum_{m=2}^k {k \choose m} x^{k-m-1} y^m 
= \frac{(x+y)^k}{x} - x^{k-1} - kx^{k-2}y,
$$
with $x = n-p$, $y = p_{c+1}$. We obtain
$$
S = p_1 \dots p_c \frac{(k+c-1)!}{k!}
\left\{
k \left[
\frac{(n-p')^{k-1}}{n-p} - (n-p)^{k-2}
\right] (n-p') - 
\right.
\qquad \qquad
$$
$$
\qquad \qquad
\left.
-k \left[
\frac{(n-p')^k}{n-p} - (n-p)^{k-1} - k p_{c+1} (n-p)^{k-2}
\right]
\right\}
$$
$$
= p_1 \dots p_c \frac{(k+c-1)!}{(k-1)!}
(n-p)^{k-2}
\left\{
-(n-p')+(n-p)+kp_{c+1}
\right\}
$$
$$
= p_1 \dots p_{c+1} \frac{(k+c-1)!}{(k-2)!} (n-p)^{k-2}.
$$

This is precisely the formula of Theorem~\ref{Thm:main}
with $c$ replaced by $c+1$. The theorem is proved.


\begin{thebibliography}{99}

\bibitem{Arnold} {\bf V. I. Arnold.}
{\em Topological classification of complex trigonometric polynomials and
the combinatorics of graphs with an identical number of vertices and
edges.} (Russian) Funktsionalnyi Analiz i Prilozheniia, {\bf 30} (1996), 
No.~1, 1--17, 96; translation in Functional Analysis and Applications, 
{\bf 30} (1996), No.~1, 1--14.

\bibitem{ELSV} {\bf T. Ekedahl, S. K. Lando, M. Shapiro,
A. Vainshtein.} {\em Hurwitz numbers and intersections
on moduli spaces of curves.} -- Inventiones Mathematicae,
{\bf 146} (2001), 297--327, {\tt arxiv:math.AG/0004096}.

\bibitem{GouJac} {\bf I. P. Goulden, D. M. Jackson.}
{\em The combinatorial relationship between trees, cacti
and certain connection coefficients for the symmetric group.}
-- European Journal of Combinatorics, {\bf 13} (1992), 357--365.

\bibitem{GouSch} {\bf A. Goupil, G. Schaeffer.}
{\em Factoring $n$-cycles and counting maps of given
genus.}
-- European Journal of Combinatorics, {\bf 19} (1998), No.~7, 
819--834.

\bibitem{Hurwitz} {\bf A. Hurwitz.}
{\em \"Uber Riemann'sche Fl\"achen mit gegebenen
Verzweigungs\-punkten.} -- Mathematische Annalen, {\bf 39} (1891),
1--61.

\bibitem{knizhka} {\bf S. K. Lando, A. K. Zvonkin.}
{\em Graphs on Surfaces and Their Applications.}
-- Springer-Verlag, 2004. 

\bibitem{LanZvo} {\bf S. K. Lando, D. Zvonkine.}
{\em On multiplicities of the Lyashko-Looijenga mapping
on the discriminant strata.}
-- Functional Analysis and its Applications,
{\bf 33} (2000), 178--188.

\bibitem{PanZvo} {\bf D. Panov, D. Zvonkine.}
{\it Counting meromorphic functions with critical points
of large multiplicities.}
-- Notes of the Scientific Seminar of the Mathematical
Department of St-Petersburg Steklov Institute (POMI),
{\bf 292} (2002), Representation Theory, Dynamical Systems,
Combinatorics, and Algorithmic Methods, No.~7, 92--119.
Available on {\tt arXiv: math.CO/0209013}.

\bibitem{Zvonkine} {\bf D. Zvonkine.}
{\em An algebra of power series arising in the intersection 
theory of moduli spaces of curves and in the enumeration of 
ramified coverings of the sphere.}
-- {\tt arXiv:math.AG/0403092} (2004).

\end{thebibliography}
\end{document}